\documentclass[11pt]{amsart}
\usepackage{amsmath,amsfonts}

\def\1ox{{ \Omega^1_{\scriptstyle{X}} }}
\def\2ox{{ \Omega^2_{\scriptstyle{X}} }}

\def\ok1{{ \Omega^1_K }}
\def\ok2{{ \Omega^2_K }}

\def\ra{{ \rightarrow }}

\def\B{{ B_{\mbox{cris}} }}

\def\Q{{ \mathbb{Q} }}

\def\8{{ {\infty } }}

\def\Gal{{ \mbox{Gal} }}

\def\^{{ ^{\wedge} }}

\newcommand{\fil}{{\rm Fil}}
\newcommand{\End}{{\rm End}}

\theoremstyle{plain}
\newtheorem{thm}{Theorem}

\numberwithin{thm}{section}

\theoremstyle{remark}
\newtheorem{remark}[thm]{Remark}

\usepackage{amsmath}
\usepackage{amsfonts}
\usepackage{amscd}
\usepackage{amssymb}

\def\Qp{\mathbb{Q}_p}

\def\Aut{{ \mbox{Aut} }}

\title{A remark on potentially semi-stable representations of
Hodge-Tate type (0,1)}

\author{Kirti Joshi and Minhyong Kim } 
\address{School of Mathematics,
Tata Institute of Fundamental Research, Homi Bhabha Road, Mumbai
400005, INDIA; Department of Mathematics, University of Arizona,
Tucson, AZ 85721, U.S.A.}  
\email{kim@math.arizona.edu, kirti@math.arizona.edu}
\begin{document}
\maketitle
\markboth{Kirti Joshi and Minhyong Kim}
{Semi-stable representations of Hodge-Tate type (0,1)} 

\section{Introduction}
The purpose of this note is to complement {\em
part} of a theorem 
from the remarkable paper of Fontaine and Mazur
on geometric Galois representations \cite{FM}.

Fix a prime $p$, and let $K$ be a finite extension of the $p$-adic numbers
$\Qp$. Fix an algebraic closure $\bar{K}$ of $K$
and let $G$ be the Galois group of $\bar{K}$ 
over $K$.
Fontaine's theory \cite{FP} classifies various types of
representations  $\rho:G\ra \mbox{Aut(V)}$ on
finite-dimensional $\Qp$-vector spaces
$V$, and we refer to op. cit. for terminology.

The theorem (C2. (ii) $\Leftrightarrow $ (iii)) 
in question from \cite{FM} says the following:

\begin{em}
If $p\geq 5$ and $(V,\rho)$ is a two-dimensional irreducible
Hodge-Tate representation of Hodge-Tate type
(0,1), then 
$\rho$ is potentially semi-stable if and only
if it is potentially crystalline.
\end{em}

	Of course, one should emphasize that the rest of the theorem
gives much more detail, namely, a complete list of possibilities, and
the theorem to follow is by no means a substitute for the refined
statements.  However, it might be worth remarking that at least this
part admits an entirely simple proof in greater generality. We note
also that \cite{FM} C2. (i) $\Leftrightarrow $(ii), the equivalence
between crystalline representations of Hodge-Tate type $(0,1)$ and
Barsotti-Tate representations, has been proved for $p\neq 2$ and
arbitrary dimension by \cite{F}, \cite{FL} in the small ramification
case and \cite{B} in general.

We are grateful to the referee for suggesting improvements.

\section{Main Theorem}
	The main result of this note is the following.
\begin{thm}\label{main}
Let $(V,\rho)$ be a irreducible finite-dimensional
Hodge-Tate representation of $G$ of Hodge-Tate type
(0,1). Then
$V$ is potentially semi-stable if and only if 
$V$ is potentially crystalline.
\end{thm}

\begin{remark}
 Note that there is no restriction on $p$ or the dimension.
\end{remark}

\begin{remark}
	It is not necessary to restrict oneself to the case of finite
residue field. It suffices to work with perfect residue fields. The
proof given below goes through without any difficulty.
\end{remark}

\begin{proof}
	Of course we need only prove the `only if' part. So assume
$\rho$ is potentially semi-stable and let $H\subset G$ be a normal
subgroup of finite index such that $\rho |H$ is semi-stable. We claim
that the $H$-representation has a non-zero crystalline
subrepresentation. To see this let $M=D_{st}(V)$ be the associated
filtered $(\phi ,N)$-module. Since $M$ is weakly-admissible, the
slopes of $\phi$ are all in the interval $[0,1]$. The equality $$p\phi
N= N\phi$$ implies that $N$ lowers the slopes by 1. Suppose there is a
non-trivial slope-zero part, $M_0$. Then $M_0$ is killed by $N$ (in
particular, it is stabilized by $N$) and hence has the structure of a
sub-module. Now $t_N(M_0)=0$ while $t_H(M_0)\leq t_N(M_0)$ from the
weak admissibility of $M$. So we necessarily have
$t_H(M_0)=0=t_N(M_0)$. Thus, $M_0$ is a weakly-admissible submodule of
$M$ on which $N=0$, so it corresponds to a crystalline
subrepresentation of $V$.  Otherwise, the slopes of $M$ are in $(0,1]$
and thus, $N$ must kill all of $M$ and $V$ is therefore crystalline.

	Let $W\subset V$ be the maximal $H$-subrepresentation of $V$
which is crystalline. This exists because the category of crystalline
representations as a sub-category of $H$-representations is closed
under direct sums and quotient objects, and hence, under taking sums
of subspaces.  The point is that $W$ is actually $G$-stable, implying
that $W=V$, and hence, that $V$ is potentially crystalline.  We see
this as follows: Let $g\in G$ and consider the subspace $g^{-1}W
\subset V$. Since $H$ is a normal subgroup, $g^{-1}W$ is an
$H$-subspace.

Claim: $g^{-1}W$ is $H$-crystalline.

	For any representation $\pi:H \ra \Aut (U)$, denote by $U_g$
the representation given by the $H$-action $\pi \circ c_g$ on the same
vector space $U$, where $c_g(h)=ghg^{-1}$. Then as an
$H$-representation, $g^{-1}W$ is isomorphic to $W_g$.

	 Thus, we need to consider the representation $W_g\otimes
\B$. But since $\B$ as an $H$-module is the restriction of a
$G$-module, $(\B)_g\simeq \B$. Therefore,
$$W_g\otimes \B \simeq W_g\otimes (\B)_g \simeq (W\otimes \B)_g$$
Since the images of the original and the
twisted actions in the automorphism group of
$W \otimes \B$ are the same, we get
$$[(W\otimes B)_g]^H=(W\otimes B)^H$$
and therefore, $W_g$ is $H$-crystalline.
Since $W$ is the maximal $H$-crystalline
subrepresentation of $V$, we must have $g^{-1}W \subset W$.
So $W$ is $G$-stable.
\end{proof}

\begin{remark}
	The idea of using kernel of $N$ also occurs in Corollary 5.3.4
of \cite{B}.
\end{remark}

\section{An example}

An important ingredient in the proof is the fact that a semi-stable
representation of Hodge-Tate type (0,1) is either crystalline or has a
{\em crystalline filtration}, i.e., a filtration whose associated
graded objects are crystalline. That this cannot be a general
phenomenon was remarked by Jannsen in \cite{J}, where he considers a
$\Q_5$ representation of dimension 2 associated to a modular form of
weight 4 and level 5.  Here we expand on a construction of \cite{FM}
to illustrate that the conditions of the theorem cannot be relaxed in
any dimension. In fact, the examples will be {\em strongly
irreducible}, that is, irreducible even when restricted to a finite
index subgroup, and a fortiori cannot have a crystalline filtration
even after such a restriction.

 The Fontaine-Mazur construction goes as follows: Let $M=<e_1,e_2>$ be
a two dimensional $\Qp$ vector space spanned by vectors $e_1,e_2$. Let
$s\geq 3$ be an odd integer and let $b\in \Qp$ be a $p$-dic integer
such that $v_p(b)=(s-1)/2$.  We define a filtration on $M$ by $\fil^i
M=M$ for $i\leq 0$, and $\fil^iM=\fil^sM=<e_1>$ for $1\leq i\leq s$;
and $\fil^iM=0$ for $i>s$.  The operator $\phi$ is defined by
$\phi(e_1)=pbe_1$ and $\phi(e_2)=be_2$; $N(e_1)=e_2$.  We see right
away that $t_H(M)=s$ and $t_N(M)=1+(s-1)/2+(s-1)/2=s$.  Now, the only
$(\phi, N)$-invariant subspace of $M$ is $<e_2>$, to which the
filtration restricts trivially. Thus,
$t_N(<e_2>)=(s-1)/2>t_H(<e_2>)=0$ and therefore, $M$ is weakly
admissible, while it does not admit any weakly admissible submodules.
Thus, $V=V_{st}(M)$ is an irreducible semi-stable representation for
$G=\Gal(\bar{\Qp}/\Qp)$ \cite{FC}.  Since $N\neq 0$, $V$ is not
crystalline, so $V$ is a two-dimensional example of a representation
which does not admit a crystalline filtration. But exactly the same
argument shows that the filtered $(\phi,N)$-module $M\otimes K$ has no
weakly admissible submodules for any extension $K$ of $\Qp$ (because
of $N$, extending the coefficients does not allow us to pick up any
other invariant lines).  Extending the coefficients of $M$ corresponds
to restricting to a subgroup of $G$, therefore, we see that $V$ is
even irreducible when restricted to finite index subgroups of $G$.

We can say more by investigating the endomorphisms of $M$. Suppose $f$
is such an endomorphism.  By checking the conditions of commuting with
$\phi$ and $N$, we see that $f$ must be a scalar $a\in \Qp$. That is,
$\End (M)=\Qp$. In fact, we can see that $\End (M\otimes K)=\Qp$ for
any extension $K$ of $\Qp$: Write $f(e_1)=xe_1+ye_2$ and
$f(e_2)=ze_1+we_2$ for $x,y,z,w \in K$. The condition of commuting
with $N$ gives us $z=0$ and $x=w$.  $f(\phi(e_2))=\phi(f(e_2))$ says
that $bw=b\sigma (w) \Rightarrow w=\sigma (w)$, where $\sigma:K
\rightarrow K$ is the Frobenius (lift) of $K$, so we get $w\in
\Qp$. Similarly, $f(\phi(e_1))=\phi(f(e_1))$ gives us
$bpwe_1+bpye_2=bpwe_1+\sigma(y)be_2$ so that $\sigma (y)=py$.  Since
$\sigma$ cannot change the valuation of non-zero elements, this
implies $y=0$. Note that the above computation does not require the
filtration at all.  By the equivalence of categories given by $V_{st}$
and $D_{st}$, we see that for any finite index subgroup $H$ of $G$,
$\End_{\Qp[H]}(V)=\Q_p$. Since we already know that $V$ is $
\Qp[H]$-simple (this {\em does} require the filtration), this implies
that $\mbox{Image}(\Qp[H])=\End_{\Qp}(V)$ by Wedderburn's theorem. In
particular, $H$ must act irreducibly on ${\rm Sym}^n(V)$ for all
$n>0$. This way, we obtain a $G$-representation for all dimensions
$n\geq 2$ which cannot admit a crystalline filtration for any finite
index subgroup.

\begin{remark}
The proof of also shows that for the semi-stable Galois module $V$
constructed above, the image of $G$ in $GL(V)$ is open. This is a
consequence of the following facts: 1) by the above calculation $V$
remains irreducible on restriction to any open subgroup of $G$ and its
endomorphism ring is $\Q_p$. 2) By a standard result the image of
Galois for any Hodge-Tate representation is open in its Zariski
closure.
\end{remark}

\vspace{3mm}

{\bf Acknowledgment:} M.K. was supported in part by a grant
from the National Science Foundation.

\end{document}